\documentclass[12pt]{amsart}
\usepackage{amssymb,latexsym,amsthm, amsmath, amsfonts, comment, bbm, fullpage, tikz-cd, wrapfig}
\usepackage[all]{xy}
\usepackage{graphicx, comment}
\usepackage[scale=1.2]{CormorantGaramond}%libertine, wrapfig}
\usepackage[utf8]{inputenc}
\usepackage[T2A]{fontenc}
\ifx\pdfoutput\undefined
\usepackage{hyperref}
\else
\usepackage[pdftex,colorlinks=true,linkcolor=blue,urlcolor=blue]{hyperref}
\fi
\theoremstyle{plain}

\theoremstyle{definition}
\newtheorem{definition}{Definition}
\newtheorem{theorem}{Theorem}
\newtheorem*{thm}{Theorem}

\newtheorem{hypo}{Hypothesis}

\newtheorem{remark}{Remark}

\newtheorem*{ack}{Acknowledgements}
\begin{document}

\title[Aperiodic and linearly repetitive Lorentz gases of finite horizon are not exponentially mixing]
      {Aperiodic and linearly repetitive Lorentz gases of finite horizon are not exponentially mixing}
      \author{Rodrigo Trevi\~no}
      \address{Department of Mathematics \\ University of Maryland \\ William E. Kirwan Hall \\ College Park, MD 20742}
      \email{rodrigo@trevino.cat}
      \author{Agnieszka Zelerowicz}
      \address{Department of Mathematics \\ University of California, Riverside\\900 University Ave.  \\ Riverside, CA 92521}
      \email{agnieszz@ucr.edu}
      \date{\today}
      \begin{abstract}
        We prove that aperiodic and linearly repetitive Lorentz gases with finite horizon are not mixing with exponential or stretched exponential speed in any dimension for any class of H\"older observables under a technical assumption known to hold in all known examples. We also bound the polynomial speed of mixing for observables in the H\"older space $H_{\alpha}$ depending on $\alpha$.
      \end{abstract}
      \maketitle
      The Lorentz gas is one of the oldest and best-studied models from statistical mechanics. It consists of a massless point particle moving through Euclidean space bouncing off a given set of scatterers $\mathcal{S}$ with ellastic collisions at the boundaries $\partial \mathcal{S}$. If the set of scatterers is periodic in space, then the quotient system, which is compact, is known as the Sinai billiard. The statistical properties of Sinai billiards in two dimensions are very well understood (see \cite{szasz:challenges, golse:survey, dettmann:survey, CM:book}) and one of the fundamental tools developed for their study are Young towers \cite{Young:annals}. Through the use of towers, Young proved that the Sinai billiard in two dimensions is mixing with exponential speed. This was an improvement of the result of Bunimovich, Sinai and Chernov \cite{BSC2}, who proved that the Sinai billiard in two dimensions was mixing with stretched exponential speed. By exponential and stretched exponential speed, we mean the following.
      \begin{definition}
        Let $F:X\rightarrow X$ be continuous map on a compact metric space $X$ which preserves a probability measure $\mu$, and let $H_\alpha$ be the set of H\"older functions on $X$ with exponent $\alpha>0$. The system $(X,F,\mu)$ is \textbf{exponentially mixing} in $H_\alpha$ if there exists a $C>0$ and $\tau\in(0,1)$ such that for any two functions $\psi_1,\psi_2\in H_\alpha$ we have
      $$\left|\int_X \psi_1\circ F^n \cdot \psi_2 \, d\mu - \mu(\psi_1)\mu(\psi_2)\right|\leq C \|\psi_1\|_\alpha \|\psi_2\|_\alpha\tau^n$$
      for all $n\geq 0$. It has \textbf{stretched exponential mixing} in $H_\alpha$ if there exist $\gamma,\tau\in(0,1)$ and $C>0$ such that for any two functions $\psi_1,\psi_2\in H_\alpha$ we have
      $$\left|\int_X \psi_1\circ F^n \cdot \psi_2 \, d\mu - \mu(\psi_1)\mu(\psi_2)\right|\leq  C \|\psi_1\|_\alpha \|\psi_2\|_\alpha\tau^{n^\gamma}$$
      for all $n\geq 0$.
      \end{definition}
      For periodic Lorentz gases in higher dimension, B\'{a}lint and T\'{o}th \cite{BT:exp} have proved essentially the only result on the speed of mixing. They prove that under certain conditions in the finite horizon case, the multi-dimensional Sinai billiard is exponentially mixing. One of the conditions, however, has not been made explicit in any construction; see \cite[\S 4]{Szasz:multi} for a more thorough discussion.
      
      In contrast to the Lorentz gas with periodic configurations of scatterers, very little is known about the behavior of the Lorentz gas with aperiodic configurations of scatterers which model quasicrystals and other low-complexity aperiodic sets. The work of Marklof-Str\"ombergsson \cite{MS:quasi, marklof:ICM} studied the Boltzmann-Grad limit, which is the limit as the size of the scatterers goes to zero, for special types of aperiodic scatterer configurations in any dimension. In \cite{TZ:ALG1} we proved a particular type of mixing for a large class of aperiodic Lorentz gases in two dimensions, but the question of speed of mixing was left open. Here we prove that exponential mixing never happens for the billiard map of aperiodic Lorentz gases with linearly repetitive scatterer configurations in any dimension assuming a technical condition that holds in all known classes of examples (see \S \ref{sec:LRS} for the definition of linear repetitivity).
\begin{figure}[t]
  \centering
  \includegraphics[width = 4.5in]{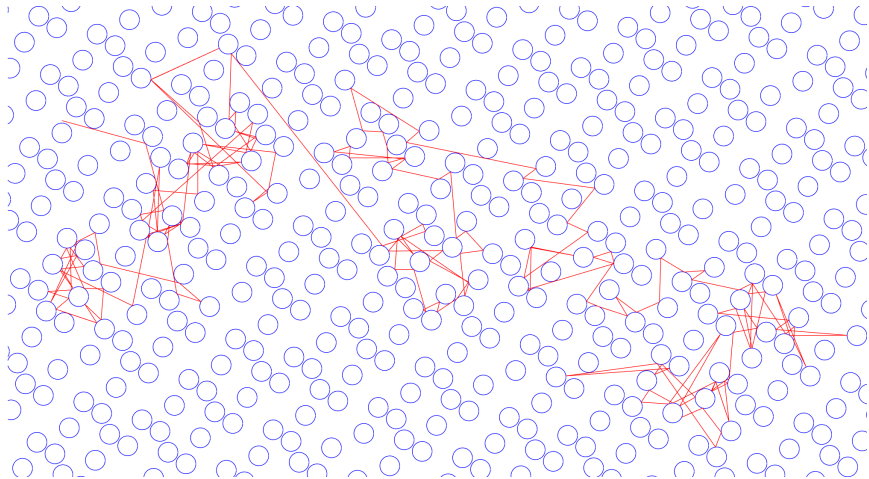}
  \caption{ A trajectory in an aperiodic Lorentz gas.}
  \label{fig:traj}
\end{figure}
      
      Let us be more specific: just as there is a compactification of a Lorentz gas with periodic scatterer configuration (which defines the Sinai billiard), there is a compactification of Lorentz gases when the scatterer configuration $\mathcal{S}$ is aperiodic and of finite local complexity. There are a couple of equivalent ways to do this (see \cite{TZ:ALG1}), but the point is that the resulting collision space $\mathcal{D}$ parametrizes not only a particular aperiodic Lorentz gas but a continuum of Lorentz gases which are locally indistinguishable to the original one. The space $\mathcal{D}$ is not a manifold but has a nice local product structure (see \S\ref{sec:ALG} below). If the billiard map $F:\mathcal{D}\rightarrow \mathcal{D}$ is mixing with respect to some measure $\mu$, and the measure $\mu$ is nice enough, then this implies the \textbf{leafwise mixing}\footnote{This type of mixing is also called \emph{global-global mixing} in \cite{lenci:infMixing, DN:mixing}.} of the billiard map $f:T^1_+(\partial\mathcal{S})\rightarrow T^1_+(\partial\mathcal{S})$ of the Lorentz gas defined on the natural Poincar\'e section $T^1_+(\partial\mathcal{S})\subset \partial\mathcal{S}\times S^{d-1}$, set of unit tangent vectors at the boundary of the scatterers which are not inwardly pointing. Denoting by $\pi:T^1_+(\partial\mathcal{S})\rightarrow \partial \mathcal{S}$ the projection to the first coordinate, leafwise mixing of the Lorentz gas in a class $\mathcal{F}$ of functions $g:T^1_+(\partial\mathcal{S})\rightarrow \mathbb{R}$ means that for any $g_1,g_2\in\mathcal{F}$:
      \begin{equation}
        \label{eqn:LeafMix}
        \lim_{n\rightarrow\infty}\lim_{R\rightarrow \infty}\frac{1}{\mathrm{Vol}(B_R(0))}\int_{\pi^{-1}(B_R(0)\cap \partial \mathcal{S})}g_1\circ f^n \cdot g_2 \, d\Xi = \beta(g_1)\beta(g_2),
      \end{equation}
      where $B_R(0)\subset \mathbb{R}^d$ is the closed ball of radius $R$ centered at $0$ and $\beta(g_i)$ is the average of $g_i$ in $\mathbb{R}^d$:
      $$\beta(g_i):=\lim_{R\rightarrow \infty}\frac{1}{\mathrm{Vol}(B_R(0))}\int_{\pi^{-1}(B_R(0)\cap \partial \mathcal{S})}g_i \, d\Xi.$$
      Here 
      %$g_i$ belong to a special class of functions which contains functions whose values depend on the local patterns of the aperiodic set underlying the aperiodic scatterer configuration $\mathcal{S}$, and
      $\Xi$ is a smooth invariant measure on the Poincar\'e section $T^1_+(\partial \mathcal{S})$ of the Lorentz gas. Exponential mixing of the billiard map $F:\mathcal{D}\rightarrow \mathcal{D}$ for functions in $\mathcal{F}$ would imply that the rate of convergence in $n$ of the outer limit above would also be exponential.
      
      In \cite{TZ:ALG1} leafwise mixing was proved for the billiard map of repetitive, aperiodic Lorentz gases of finite local complexity with finite horizon in two dimensions (in fact $K$-mixing was proved). The mixing rate was left as an open question. The main result in this paper shows that the rate of mixing cannot be exponential, and in fact it is significantly restricted by the geometry of the collision space $\mathcal{D}$ (assuming technical Hypothesis \ref{hyp:1}, which we will define at the end of Section \ref{sec:LRS}, after establishing more structure of our space; see also Remark (ii) below).   In the following, $\mu$ is the ``natural'' $F$-invariant probability measure on $\mathcal{D}$ (by that we mean absolutely continuous in the sense of Definition \ref{def:abscont}), and $H_\alpha^\perp$ is the space of transversally $\alpha$-H\"older functions on the collision space $\mathcal{D}$ defined in \S \ref{sec:Far}.
      \begin{theorem}
        The billiard map of aperiodic Lorentz gases in $\mathbb{R}^d$ with finite horizon and linearly repetitive scatterer configurations satisfying Hypothesis \ref{hyp:1} is not mixing with exponential or stretched exponential speed in any $H_\alpha^\perp$, $\alpha>0$. Moreover, if for some $\alpha>0$ and $C>0$ we have that
        $$\left|\int_\mathcal{D} \psi_1\circ F^n \cdot \psi_2 \, d\mu - \mu(\psi_1)\mu(\psi_2)\right|\leq C \|\psi_1\|_\alpha \|\psi_2\|_\alpha n^{-\gamma}$$
        for all $n\in\mathbb{N}$, $\psi_1,\psi_2\in H_\alpha^\perp$, then $\gamma\leq 2(d+\alpha)$.
      \end{theorem}
      \begin{remark}
        Some remarks:
        \begin{enumerate}
        \item This should be contrasted with the exponential mixing of Sinai billiards, revealing a strong qualitative difference between the periodic and the aperiodic cases in two dimensions.
        \item The class of aperiodic scatterer configurations for which we prove this, the so-called linearly repetitive sets, is a very broad class. In particular, it includes all scatterer configurations obtained from substitution tilings and mixed/S-adic/globally-random substitution tilings (e.g. \cite{ST:random}). These types of tilings give linearly repetitive Delone multisets which satisfy Hypothesis \ref{hyp:1}. In fact, we do not know of any linearly repetitive Delone multiset which does not satisfy this hypothesis.
        \item Linearly repetitive cut-and-project sets with cubical windows were characterized in \cite{HKW:CAP} and thus our results hold exactly for the class of sets characterized in \cite{HKW:CAP}. A further characterization of repetitivity of cubical cut-and-project sets was done in \cite{HJKW:stat}. It seems like our proof is likely to go through even for non-linearly repetitive aperiodic scatterers coming from typical cut-and-project with cubical windows using the tools from \cite{HJKW:stat}. Recently, Walton established a complete characterization of linearly repetitive cut and project sets with polytopal windows \cite{walton}.% seems to be the crucial estimate we need in our proof here, which is (\ref{eqn:balance}) here in the case of linearly repetitive Delone sets.
        \item The collision spaces for which these results hold are not locally Euclidean (manifolds). As such, non-constant functions with H\"older exponents $\alpha>1$ exist, and this exponent in spirit describes the regularity in the non-Euclidean direction of the space. The bound on $\gamma$ depending on $\alpha$ relies on the metric used to obtain this bound in terms of the H\"older exponent. This bound is for the natural tiling metric typically used in spaces of repetitive tilings of finite local complexity.
        \end{enumerate}
      \end{remark}
      So what does slow mixing for functions which are transversally $\alpha$-H\"older look like from the point of view of the Lorentz gas on $\mathbb{R}^d$? If $g\in H_\alpha^\perp$ then there exists $\hat{g}:T^1_+(\partial\mathcal{S})\rightarrow \mathbb{R}$ and 
$C_{\hat{g}}= |g|^\perp_{\alpha}\geq0$ such that if, for two translation equivalent scatterers $S_1, S_2\subset \mathcal{S}$ with $\varphi_t(S_1)=S_2$, the scatterer configurations inside a ball of radius $R$ around them are translation equivalent, then $|\hat{g}(x_1)-\hat{g}(x_2)|\leq C_{\hat{g}} R^{-\alpha}$ for any $(x_1,x_2)\in \pi^{-1}(\partial S_1)\times\pi^{-1}(\partial S_2)$ satisfying $\varphi_t(x_1)=x_2$. For two functions $\hat{g}_1,\hat{g}_2:T^1_+(\partial\mathcal{S})\rightarrow \mathbb{R}$ of this type, it follows from \cite{TZ:ALG1} that they satisfy the leafwise mixing (\ref{eqn:LeafMix}). Our main result here says that if there are $C,\gamma>0$ such that
      \begin{equation*}
        \begin{split}
          &\left|\lim_{R\rightarrow \infty}\frac{1}{\mathrm{Vol}(B_R(0))}\int_{\pi^{-1}(B_R(0)\cap \partial \mathcal{S})}\hat{g}_1\circ f^n \cdot \hat{g}_2 \, d\Xi - \beta(\hat{g}_1)\beta(\hat{g}_2)\right| \\ &\hspace{3in}\leq C (\|\hat{g}_1\|_{C^0}+C_{\hat{g}_1}) (\|\hat{g}_2\|_{C^0}+C_{\hat{g}_2})n^{-\gamma}
        \end{split}
      \end{equation*}
      for all $n>0$, then $\gamma\leq 2(d+\alpha)$.
      
      This paper is organized as follows. In \S \ref{sec:backDelone} we review aperiodic Delone multisets and the structure of pattern spaces. In \S \ref{sec:ALG} we recall the setup of aperiodic Lorentz gases from \cite{TZ:ALG1} and their relation to aperiodic Delone multisets. In \S \ref{sec:LRS} we recall relevant facts about the structure of linearly repetitive Delone multisets, following mostly the results of \cite{AC:tower, ACCDP:LR}, and state our technical Hypothesis \ref{hyp:1}. Sections \ref{sec:Far} and \ref{sec:nonFast} show how the structure of pattern spaces of linearly repetitive aperiodic sets prohibits fast mixing, that is, they give the proof of the main result.
      \begin{ack}
        We thank Daniel Coronel for helpful insights on our Hypothesis \ref{hyp:1}. This work was partially supported by the Simons Collaboration Grant 712227.
      \end{ack}
      \section{Background: Delone multisets and repetitivity}
      \label{sec:backDelone}
      A \textbf{Delone set} in $\mathbb{R}^d$ is a countable subset $\Lambda\subset \mathbb{R}^d$ that is
      \begin{description}
      \item[uniformly discrete] there exists an $r_\Lambda>0$ such that, for any $x\in \mathbb{R}^d$, $B_{r_\Lambda}(x)\cap \Lambda$ contains at most one point, and
      \item[relatively dense] there exists an $R_\Lambda>0$ such that, for any $x\in\mathbb{R}^d$, $B_{R_\Lambda}(x)\cap \Lambda$ contains at least one point.
      \end{description}
      These numbers are called the \textbf{packing} and \textbf{covering} radii, respectively. A \textbf{Delone multiset} is a set $\Lambda$ of the form 
      $$\Lambda = \bigcup_{i\in \mathfrak{C}} \Lambda_i\times \{i\}\subset \mathbb{R}^d\times \mathfrak{C},$$
      where each $\Lambda_i$ is a Delone set, such that the Delone sets $\Lambda_i$ are pairwise disjoint, and $\mathfrak{C}$ is a finite set. In this case the index $i$ of a set $\Lambda_i$ can be thought of as a label or color. The translation of a Delone multiset $\Lambda$ by $t\in\mathbb{R}^d$ is denoted by
      $$\varphi_t(\Lambda) := \bigcup_{i\in\mathfrak{C}}\left(\Lambda_i-t\right)\times\{i\}.$$
      The multiset $\Lambda$ is \textbf{aperiodic} if $\varphi_t(\Lambda)=\Lambda$ implies that $t = 0$. Denote by $\pi: \Lambda\rightarrow \mathbb{R}^d$ the projection of $\Lambda$ to the first coordinate. $\pi(\Lambda)$ is a Delone set.
      
      A \textbf{cluster} of $\Lambda$ is a finite subset $C\subset \Lambda$. An \textbf{$R$-cluster} is a cluster of the form $C = \Lambda\cap B_R^*(p)$ for some $p\in\Lambda_j$ for some $j$, where $B^*_R(p) := B_R(p) \times \mathfrak{C}$. A Delone multiset $\Lambda$ has \textbf{finite local complexity} (FLC) if for any $R>0$ the collection of all $R$-clusters, up to translation equivalence, is finite. A Delone multiset $\Lambda$ is \textbf{repetitive} if for any $R>0$ there exists a $T>0$ such that for any $R$-cluster $C\subset \Lambda$ and $p\in\Lambda_j$ for some $j$, the $T$-cluster around $p$ contains point $p'$ such that $B_R^*(p')\cap \Lambda$ is translation equivalent to $C$. 

      For a Delone multiset $\Lambda$ and some translates $\varphi_s(\Lambda), \varphi_t(\Lambda)$ consider the quantity
      \begin{equation*}
        \begin{split}
          &d^*(\varphi_s(\Lambda),\varphi_t(\Lambda)) \\
          &\hspace{.5in}:= \inf_{\varepsilon>0} \{\mbox{ there exist $x,y\in B_\varepsilon(0)$ such that }B_{1/\varepsilon}^*(0)\cap \varphi_{s+x}(\Lambda) = B_{1/\varepsilon}^*(0)\cap \varphi_{t+y}(\Lambda)\},
        \end{split}
      \end{equation*}
      and introduce the function $d$ on the set of pairs of translates
      \begin{equation}
        \label{eqn:distance}
        d(\varphi_s(\Lambda), \varphi_t(\Lambda)) = \min\{d^*(\varphi_s(\Lambda),\varphi_t(\Lambda)), 2^{-1/2}\}.
      \end{equation}
      This function is a metric on this set \cite{LMS:PP}, and we consider the completion, with respect to this metric, of the set of all translates of $\Lambda$:
      $$\Omega_\Lambda := \overline{\left\{\varphi_t(\Lambda):t\in\mathbb{R}^d\right\}}.$$
      This set is called the \textbf{pattern space of $\Lambda$}. If $\Lambda$ has finite local complexity, which will be assumed throughout the paper, then $\Omega_\Lambda$ is a compact metric space with an $\mathbb{R}^d$ action defined by the translation $\varphi_t$. If $\Lambda$ is repetitive -- which implies FLC -- then the action is minimal.

      Let $\Lambda$ be a repetitive Delone multiset of finite local complexity. For $\Lambda'\in\Omega_\Lambda$ and $R>0$, let $C_{\Lambda',R}$ be the \textbf{cylinder} set defined by the $R$-neighborhood of the origin in $\Lambda'$, that is,
      $$C_{\Lambda',R}:= \{\Lambda''\in\Omega_\Lambda: \Lambda''\cap B_R^*(0)= \Lambda'\cap B_R^*(0) \}.$$
      Any clopen subset of a cylinder in $\Omega_\Lambda$ is called a \textbf{local transversal}.
      
      The pattern space has a \textbf{canonical transversal} defined by the set
      $$\mho_\Lambda := \left\{\Lambda' \in \Omega_\Lambda: 0\in \pi(\Lambda') \right\}.$$
      For any $R>0$, the canonical transversal is partitioned into finitely many cylinder subsets $\mathcal{P}_R = \{\mathcal{C}_P\}$, where $P$ is an $R$-cluster, as follows. By finite local complexity, given $R>0$, there exist finitely many $R$-clusters $C_{R}^1,\dots, C_{R}^{k_R}$. As such, if $\Lambda'\in \mho_\Lambda$, then the $R$-cluster of $\Lambda'$ around the origin is exactly one of the $C_{R}^i$, and in that case $\Lambda'$ belongs to the element of the partition of $\mho_\Lambda$ determined by $C_R^i$. Thus the partition is
      $$\mathcal{P}_R = \bigsqcup_{i=1}^{k_R} \left\{\Lambda'\in\mho_\Lambda: \Lambda'\cap B_R^*(0) = C_R^i \right\}.$$
      It follows that the space $\Omega_\Lambda$ has a local product structure of the form $B_\varepsilon(0)\times \mathcal{C}$, where $\mathcal{C}$ is a Cantor set. 

            Our main result relies on transversally $\alpha$-H\"older functions, which we now define. These are continuous functions which are $\alpha$-H\"older but only in the transverse direction. Given a continuous function $f:\Omega_\Lambda\rightarrow \mathbb{R}$, define the transversally $\alpha$-H\"older seminorm for $f$ as
      \begin{equation}
        \label{eqn:HoldConst}
        \begin{split}
          |f|_\alpha^\perp&:= \sup_{\Lambda'\in\Omega_\Lambda}\sup_{\substack{\Lambda_1\neq \Lambda_2 \\ \in C_{\Lambda',R_{\Lambda}}}}\frac{|f(\Lambda_1)-f(\Lambda_2)|}{d(\Lambda_1,\Lambda_2)^\alpha}\\
          &=  \sup_{\Lambda'\in\Omega_\Lambda}\sup_{R>R_{\Lambda}}\sup_{\substack{\Lambda_1\neq \Lambda_2 \\ \in C_{\Lambda',R}}}\frac{|f(\Lambda_1)-f(\Lambda_2)|}{R^{-\alpha}}.
        \end{split}
      \end{equation}
      Define the space of transversally $\alpha$-H\"older functions as
      $$H_\alpha^\perp(\Lambda):= \{f:\Omega_\Lambda\rightarrow \mathbb{R}: f\mbox{ is continuous and }|f|^\perp_\alpha<\infty\}$$
      which, when endowed with the norm $\|f\|_\alpha = \|f\|_{C^0}+|f|^\perp_\alpha$, is a Banach space. We note that since the H\"older regularity is being controled by $\alpha$ on the local transversals, the spaces $H^\perp_\alpha$ are non-trivial for all values of $\alpha>0$, unlike H\"older spaces for manifolds where spaces of $\alpha$-H\"older functions become trivial (one-dimensional) once $\alpha$ is large enough.
      
      %An invariant measure $\mu$ for the $\mathbb{R}^d$ action on $\Omega_\Lambda$ has a local product structure $\mathrm{Leb}\times \nu$, where $\mathrm{Leb}$ is the Lebesgue measure and

      The space $\Omega_\Lambda$ carries a natural transverse $\mathbb{R}^d$-invariant measure $\nu$ in the sense of Bowen-Marcus \cite{BW:UE}. The measure $\nu$ is defined as follows. Let $\Lambda'\in\Omega_\Lambda$ and $A\subset \mathbb{R}^d\times \mathfrak{C}$ a bounded subset. For a cluster $P$, define
      $$L_P^{\Lambda'}(A):= |\{t\in\mathbb{R}^d:\varphi_t(P)\subset A\cap \Lambda'\}|$$
      to be the number of translates of $P$ contained in $A$. Taking a limit we get the frequency of a cluster
      $$\mathrm{freq}_{\Lambda'}(P) = \lim_{R\rightarrow \infty}\frac{L_P^{\Lambda'}(B_R(0))}{\mathrm{Vol}(B_R(0))}.$$
      This number assigned to $P$ is clearly $\mathbb{R}^d$ invariant in that $\mathrm{freq}_{\Lambda'}(P) = \mathrm{freq}_{\varphi_\tau(\Lambda')}(P)$ for any $\tau\in\mathbb{R}^d$. The frequency map induces a transverse $\mathbb{R}^d$-invariant measure $\nu_{\Lambda'}$ on $\mho_\Lambda$, defined for a cylinder set $\mathcal{C}_P$ of an $R$-cluster $P$,
      \begin{equation}
        \label{eqn:freqMeas}
        \nu_{\Lambda'}(\mathcal{C}_P) := \mathrm{freq}_{\Lambda'}(P).
      \end{equation}
      We call such a measure $\nu_{\Lambda'}$ a \textbf{frequency measure}. Since the frequency map is invariant under translations $\nu$ is an $\mathbb{R}^d$-invariant transverse measure. When the frequencies are independent of the set $\Lambda'$ used to compute them, the system has \textbf{uniform cluster frequencies}.

      Using the local product structure of $\Omega_\Lambda$, any $\mathbb{R}^d$ invariant transverse measure $\nu$ can be locally paired with the Lebesgue measure and obtain a $\mathbb{R}^d$-invariant measure on $\Omega_\Lambda$. In the case of uniform cluster frequencies this measure is unique and the system is uniquely ergodic \cite[Theorem 2.6]{LMS:PP}. This will be the only relevant case in this paper.
      \section{Aperiodic Lorentz Gases}
      \label{sec:ALG}
      Let $\mathcal{S}\subset \mathbb{R}^d$ be a countable collection of open convex topological balls with $C^{d-1}$ boundaries having pairwise disjoint closures, called a set of \textbf{scatterers}. The \textbf{Lorentz gas} on $\mathbb{R}^d$ with scatterers $\mathcal{S}$ is the system describing the free motion of a point particle in $\mathbb{R}^d$ having elastic collisions at $\partial \mathcal{S}$. That is, it is the flow on $\mathbb{R}^d\times T^1\mathbb{R}^d$ defined by geodesics on $(\mathbb{R}^d\setminus \mathcal{S}) \times T^1\mathbb{R}^d$ and changing directions at $\partial\mathcal{S}$. The system has \textbf{finite horizon} if there is a constant $M$ such that the time between any two collisions is bounded by $M$. In this case, the set
      $$T^1_+(\partial\mathcal{S}):= \left\{( x,v)\in\partial \mathcal{S}\times\mathbb{R}^d: \|v\|=1 \mbox{ and there is $\epsilon>0$ such that } x+tv\not\in \mathcal{S}\mbox{ for }t\in (0,\epsilon)\right\}$$
      serves as a Poincar\'e section for this flow and the Poincar\'e map for this section is called the \textbf{billiard map} $f:T^1_+(\partial\mathcal{S})\rightarrow T^1_+(\partial \mathcal{S})$.

      To any collection of scatterers we will assign a countable set $\Lambda_\mathcal{S}$ as follows. Define an equivalence relation on the connected components of $\mathcal{S}$ by $S_i\sim S_j$ if they are translation equivalent, that is, if $\varphi_s(S_i) = S_j$ for some $s\in\mathbb{R}^d$. Let $\mathfrak{C}$ be the set of classes of translation-equivalent scatterers. For each $i\in \mathfrak{C}$, let $\Lambda_i^\mathcal{S}$ be the union of centers of masses of components in the class $i$, and define
      $$\Lambda_\mathcal{S} := \bigsqcup_{i\in\mathfrak{C}}\Lambda_i^\mathcal{S}.$$
      \begin{definition}
        The scatterer configuration $\mathcal{S}$ is aperiodic, repetitive, has finite local complexity, etc, if the associated set $\Lambda_\mathcal{S}$ is a Delone multiset with the same properties.
      \end{definition}

\begin{remark}
It can be observed that a Delone multiset $\Lambda_\mathcal{S}$ being repetitive and of finite local complexity implies those properties for each Delone set $\Lambda_i^\mathcal{S}$. On the other hand, $\Lambda_\mathcal{S}$ may be aperiodic even if some of the sets
$\Lambda_i^\mathcal{S}$ are periodic.
\end{remark}

      \begin{remark}
        It should be emphasized that whenever a scatterer configuration $\mathcal{S}$ is aperiodic and of finite local complexity then the associated set $\Lambda_\mathcal{S}$ is a Delone multiset consisting of \textbf{finitely many classes} parametrized by a \text{finite} set $\mathfrak{C}$.
      \end{remark}
      Suppose that $\mathcal{S}$ is aperiodic, repetitive, and of finite local complexity. Then we can identify it with a subset $\Sigma'\subset \Omega_\mathcal{S} := \Omega_{\Lambda_\mathcal{S}}$ as follows. 
%Define $R_\mathcal{S} = \{t\in\mathbb{R}^d:0\in\varphi_{t}(\mathcal{S})\}$ to be the set of translates of $\mathcal{S}$ which cover the origin and note that this contains $\Lambda_\mathcal{S}$. This defines a subset $\Sigma' = \{\varphi_t(\Lambda_\mathcal{S}):t\in R_\mathcal{S}\}\subset \Omega_{\mathcal{S}}$. 
%\az{$\mathcal{S}$ can be identified with $\Sigma' $ (a single Lorentz gas), while $\Sigma_\mathcal{S}$ is the whole collision space (continuum of Lorentz gasses).}
Observe that $t\in\mathbb{R}^d$ is in $\mathcal{S}$ if and only if $\varphi_t(\mathcal{S})$ contains $0$. As such, $\mathcal{S}$ can be identified with a set of vectors $R_\mathcal{S} = \{t\in\mathbb{R}^d:0\in\varphi_{t}(\mathcal{S})\}$,
which in turn can be identified with the set of translates of $\mathcal{S}$ that cover the origin or with the following subset of translates of $\Lambda_\mathcal{S}$, $$\Sigma' = \{\varphi_t(\Lambda_\mathcal{S}):t\in R_\mathcal{S}\}\subset \Omega_{\mathcal{S}}.$$

The above set is defined along the orbit of $\Lambda_\mathcal{S}$, so taking the closure in $\Omega_{\mathcal{S}}$ we obtain the set $\Sigma_\mathcal{S} = \overline{\Sigma'}$. As such, the Lorentz gas can be represented as a flow $\Phi_t:T^1_+(\overline{\Omega_\mathcal{S}\setminus \Sigma_\mathcal{S}})\rightarrow T^1_+(\overline{\Omega_\mathcal{S}\setminus \Sigma_\mathcal{S}})$, where
      $$T^1_+(\overline{\Omega_\mathcal{S}\setminus \Sigma_\mathcal{S}}):= \left(\overline{\Omega_\mathcal{S}\setminus \Sigma_\mathcal{S}}\times T^1\mathbb{R}^d\right)/\sim,$$
      where, for $\Lambda\in\partial \Sigma_{\mathcal{S}}$, $(\Lambda,v^-)\sim (\Lambda,v^+)$, where $v^-$ is the incoming vector in the collision determined by $\Lambda\in\partial\Sigma_{\mathcal{S}}$, and $v^+$ is the unique outgoing vector determined by ellastic collision. The flow has the property that, if $\Lambda\not\in \partial \Sigma_\mathcal{S}$ and $(\Lambda,v) \in \overline{\Omega\setminus \Sigma_\mathcal{S}}\times T^1\mathbb{R}^d $, then
      $$\Phi_t((\Lambda, v)) = (\varphi_{vt}(\Lambda), v)$$
      for $t$ in a maximal interval $(t_-,t_+)$ such that $\varphi_{vt}(\Lambda)\cap \partial\Sigma_\mathcal{S}= \varnothing$. Then the definition of $T^1_+(\overline{\Omega_\mathcal{S}\setminus \Sigma_\mathcal{S}})$ determines how to continue a trajectory once $\varphi_{vt_+}(\Lambda)\in \partial\Sigma_\mathcal{S}$.

      The flow $\Phi_t$ is thus defined to be the unique flow so that the aperiodic Lorentz gas evolving as a system on 
%$(\mathbb{R}^d\times \mathcal{S})\times T^1\mathbb{R}^d$ 
$(\overline{\mathbb{R}^d\setminus \mathcal{S}})\times T^1\mathbb{R}^d$ 
can be ``seen'' as a flow evolving inside of $\overline{\Omega_\Lambda\setminus \Sigma_\mathcal{S}}\times T^1\mathbb{R}^d$. In this space, which involves the translation closure of $\Lambda_\mathcal{S}$, not only do we find the Lorentz gas defined by $\mathcal{S}$ but a continuum of Lorentz gases defined on scatterer configurations which are locally indistinguishable from $\mathcal{S}$.

      The local product structure of $\Omega_\mathcal{S}$ gives us a nice system of coordinates for the Poincar\'e map:
      \begin{equation}
        \label{eqn:coordinates}
        \mathcal{D}:= \bigsqcup_{i\in\mathfrak{C}}T^1_+\left(\partial S_i \right) \times \mathcal{C}_i 
      \end{equation}
      where
      \begin{equation}
        \label{eqn:goodTangent}
        T^1_+\left(\partial S_i \right):= \left\{( x,v)\in\partial S_i\times\mathbb{R}^d:\|v\|=1 \mbox{ and there is $\epsilon>0$ such that } x+tv\not\in \mathcal{S}\mbox{ for }t\in (0,\epsilon)\right\}.
      \end{equation}
      $S_i$ is a representative of the scatterer class $i$ and $\mathcal{C}_i\subset \mho_\mathcal{S}$ is the subset of the canonical transversal having a point of $\Lambda_i^\mathcal{S}$ at the origin. Thus a point in $\mathcal{D}$ records what type of scatterer it is on (the index $i$), where the center of mass of this scatterer is located in $\Lambda_{\mathcal{S}}$ (the coordinate $ \mathcal{C}_i$), where on the scatterer it is (the coordinate in $\partial S_i$), and in which direction it is going (the vector $v$). The set $\mathcal{D}$ is then the \textbf{collision space} and the Poincar\'e map $F:\mathcal{D}\rightarrow \mathcal{D}$ is defined on this set.

\begin{definition}\label{def:abscont}
An invariant measure $\mu$ for $F$ is \textbf{absolutely continuous} if it is locally absolutely continuous with respect to $L_1\times \nu\times L_2$, where $L_1$ is a Lebesgue measure on $\partial S_i$, $\nu$ is a frequency measure defined by (\ref{eqn:freqMeas}), and $L_2$ is a Lebesgue measure on $S^{d-1}_+$.
\end{definition}

\begin{remark}
In \cite{TZ:ALG1} we proved that any measure satisfying above definition has $K$ property and in particular is ergodic. In the current paper we in addition assume  unique ergodicity of $\mathbb{R}^d$-action on $\Omega_{\Lambda}$.In this case, and because any two ergodic probability measures are singular, there is a unique pobability measure $\mu$ satisfying the above definition.
\end{remark}

      \section{The structure of linearly repetitive sets}
      \label{sec:LRS}
      In this section we recall some properties of linearly repetitive Delone sets which will be crucial to the proof of the main result. Our main reference for this section is \cite{AC:tower}. Note that \cite{AC:tower} is written entirely in terms of Delone sets, whereas here we will need analogous results for Delone multisets. The extension for Delone multisets is rather straight forward; we will comment on the extensions at every step.

      The \textbf{repetitivity function} $M_\Lambda$ assigns to each $R>0$ the least $T$ such that every ball of radius $T$ contains a copy of every $R$-cluster that is found in $\Lambda$. A Delone set (or Delone multiset) $\Lambda$ is \textbf{linearly repetitive} if there 
exists an $L>1$ such that the repetitivity function is bounded as $M_\Lambda(R)\leq L R$. It is known that linear repetitivity implies unique ergodicity \cite{ACCDP:LR} of the $\mathbb{R}^d$ action on $\Omega_\Lambda$ but this can fail for very mild superlinear repetitivity \cite{LP:rep, CN:examples}.
%      \begin{lemma}[See \cite{ACCDP:LR}]
%        \label{lem:rR}
%        Let $\Lambda$ be a linearly repetitive aperiodic Delone set with constant $L>1$. Then for any cluster of the form $P = B_R(p)\cap \Lambda$ with $p\in \Lambda$ and $R>0$,
%        $$\frac{R}{L+1}\leq r_{P}\leq R_P \leq LR.$$
%      \end{lemma}
%      \begin{corollary}
%        \label{cor:freqBnds}
%        There exist a constant $C>1$ such that for any cluster of the form $P = B_R(p)\cap \Lambda$ if $\mathcal{C}_P$ is the transversal cylinder set associated to the cluster $P$, we have that
%        $$C^{-1}R^{-1}\leq \nu(\mathcal{C}_P)= \mathrm{freq}(P)\leq CR^{-1}.$$
%      \end{corollary}

 For any local transversal $C$, its \textbf{recognition radius} is defined as
      $$\mathrm{rec}(C):= \inf \{R>0: C_{\Lambda',R}\subset C\mbox{ for all $\Lambda'\in C$}\}$$
      which is finite for any repetitive Delone set (or multiset) of finite local complexity. If $C$ is a local transversal and $D\subset\mathbb{R}^d$ an open set, define
      $$C[D]:= \{\varphi_t(\Lambda'): t\in D,\;\; \Lambda'\in C\}.$$
      This type of set is called a \textbf{box} if the map from $D\times C$ to $B = C[D]$ given by $(t,\Lambda')\mapsto \varphi_t(\Lambda')$ is a homeomorphism. If $C$ is a local transversal, then there exists $r(C)>0$ such that $C[D]$ is a box provided that $D\subset B_{r(C)}(0)$. A \textbf{box decomposition} of $\Omega_\Lambda$ is a finite collection of pairwise disjoint boxes $\mathcal{B} = \{B_1,\dots, B_k\}$ such that the closures of the boxes form a cover of $\Omega_\Lambda$. In this case, the boxes will be denoted by $C_i[D_i]$, where $C_i$, called the \textbf{base}, is contained in $B_i$ and $D_i$ contains the origin.

      \begin{definition}
        A box decomposition $\mathcal{B}' = \{C'_i[D_i']\}_{i=1}^{k'}$ is \textbf{zoomed out} of another box decomposition $\mathcal{B} = \{C_i[D_i]\}_{i=1}^{k}$ if the following properties hold between them:
        \begin{enumerate}
        \item If $\Lambda'\in C_i'$ satisfies $\varphi_{x'}(\Lambda')\in \varphi_y(C_j)$ for some $x'\in \overline{D_i'}$ and $y\in \overline{D_j}$, then $\varphi_{x'}(C_i')\subset \varphi_y(C_j)$;
        \item if $x'\in \partial D_i'$, there exist $j$ and $y\in \partial D_j$ such that $\varphi_{x'}(C_i')\subset \varphi_y(C_j)$;
        \item For every box $B'\in \mathcal{B}'$ there is $B\in\mathcal{B}$ such that $B\cap B'\neq \varnothing$ and $\partial B\cap \partial B' = \varnothing$;
        \item If for each $1\leq i\leq k'$ and $1\leq j \leq k$ we define
          $$O_{i,j} = \{x'\in D_i' : \varphi_{x'}(C_i')\subset C_j \},$$
          then for each $1\leq i \leq k'$ we have that
          $$\overline{D_i'} = \bigcup_{j=1}^k\bigcup_{x'\in O_{i,j}}\varphi_{x'}(\overline{D_j})$$
          where all the sets on the right have pairwise disjoint interiors;
        \item The base of $\mathcal{B}'$ is included in the base of $\mathcal{B}$, that is, $\bigcup_i C_i'\subset \bigcup_j C_j$.
        \end{enumerate}
        It follows from these properties that we must also have that
        \begin{equation}
          \label{eqn:transvPart}
          C_j = \bigcup_{i=1}^{k'}\bigcup_{x'\in O_{i,j}}\varphi_{x'}(C_i')
        \end{equation}
        for $1\leq j \leq k$ \cite[Lemma 3.2]{AC:tower}. A \textbf{tower system} is a sequence of box decompositions $\{\mathcal{B}_n\}$ such that $\mathcal{B}_{n+1}$ is zoomed out of $\mathcal{B}_n$.
      \end{definition}
      Tower systems give Delone sets a hierarchical structure and they always exist for pattern spaces of aperiodic, repetitive Delone sets \cite{BBG}. Tower systems are a manifestation of Rokhlin-Kakutani towers for these minimal systems. As such, if $\Lambda$ is a repetitive Delone multiset of finite local complexity, it also admits a tower system. In the linearly repetitive case, they have a particular nice structure which we will use.

      If $\{C_n\}$ is a decreasing sequence of local transversals with diameter going to zero and $\{\mathcal{B}_n\}$ is a tower system, then \textbf{$\{\mathcal{B}_n\}$ is adapted to $\{C_n\}$} if we have that $\mathcal{B}_n = \{C_{n,i}[D_{n,i}]\}_{i=1}^{k_n}$ such that $C_n = \bigcup_i C_{n,i}$. In such case define
      $$O_{i,j}^{(n)} = \{x\in D_{n,i}: \varphi_x(C_{n,i})\subset C_{n-1,j}\}\hspace{.7in}\mbox{ and } \hspace{.7in} m_{i,j}^{(n)} = |O_{i,j}^{(n)}|.$$
      The numbers $m_{i,j}^{(n)}$ are the entries of a $k_{n}\times k_{n-1}$ matrix, called the \textbf{transition matrix}. Since each local transversal $C_{n,i}$ is defined by the translation-equivalence class of some cluster $P_{n,i}$ (not necessarily an $R$-cluster) with a choice of a point of the cluster being the origin, these numbers record the number of clusters of type $P_{n-1,j}$ which are found in the cluster of type $P_{n,i}$.

      Given a box decomposition $\mathcal{B} = \{C_i[D_i]\}_{i=1}^k$ its \textbf{internal and external radii} are defined respectively as
      \begin{equation*}
        \begin{split}
          r_{int}(\mathcal{B}) &= \min_{i\in \{1,\dots, k\}}\sup \{R>0:  B_R(0) \subset D_i \}, \\
          R_{ext}(\mathcal{B}) &= \max_{i\in \{1,\dots, k\}}\inf \{R>0: D_i \subset B_R(0)\}, \\
          \mathrm{rec}(\mathcal{B}) &= \max_{i\in \{1,\dots, k\}}\mathrm{rec}(C_i).
        \end{split}
      \end{equation*}
      We can now state the main theorem on the structure of linearly repetitive Delone sets.
      \begin{thm}[\cite{AC:tower}]
        Let $\Lambda$ be an aperiodic linearly repetitive Delone set with repetitivity constant $L>1$ and $0\in \Lambda$. Define $\lambda=6L(L+1)^2$, $C_n:= C_{\Lambda, \lambda^n}$ for $n>0$, and $C_0=  \mho_\Lambda  $. Then there exists a tower system 
$\{\mathcal{B}_n\}= \{C_{n,i}[D_{n,i}]\}_{i=1}^{k_n}$ of $\Omega_\Lambda$ adapted to $\{C_n\} $ such that
        \begin{enumerate}
        \item $C_{n+1}\subset C_{n,1}$ 
        \item setting
          $$K_1:= \frac{1}{2(L+1)}- \frac{L}{\lambda-1}\mbox{ and }K_2:= \frac{\lambda L}{\lambda-1}$$
          we have that $0<K_1<1<K_2$ and moreover
          $$K_1 \lambda^n\leq r_{int}(\mathcal{B}_n)< R_{ext}(\mathcal{B}_n) \leq K_2 \lambda^n,$$
        \item $\mathrm{rec}(\mathcal{B}_n)\leq (2L+1)\lambda^n$
        \end{enumerate}
        for all $n$. In addition, the transition matrices have strictly positive entries and are uniformly bounded in size and norm.
      \end{thm}
      If $\mathcal{B} = \{C_{n,i}[D_{n,i}]\}$ is a tower system adapted to some $\{C_n\} = \{C_{\Lambda,\lambda^n}\}$ as in the theorem above, we will refer to it as a \textbf{nice tower system}.
      \begin{remark}
        \label{rem:general}
        This theorem above was proved for Delone sets in \cite{AC:tower}; here we outline the necessary modifications to get the result for Delone \emph{multisets}.
        
      The proof involves the following steps (see \cite[\S 4]{AC:tower} for the precise details): first, start with the linearly repetitive set $\Lambda'\in\mho_\Lambda\subset\Omega_\Lambda$ and consider the Voronoi tiling $\mathcal{T}_0$ defined by it. This gives the first collection of boxes $\mathcal{B}_0$: 
\begin{itemize}
\item$\{D_{0,i}\}_i$ are isometric to the interiors of different tile types with labels accounted for;
\item for each $i$, $C_{0,i}\subset \mho_\Lambda$ is the clopen subset corresponding to Delone sets containing origin and such that the Voronoi cell containing the origin is $D_{0,i}$;
\item $C_{0,1}\subset \mho_\Lambda$ is the clopen subset corresponding to Delone sets containing origin and such that the Voronoi cell containing the origin is the same as for $\Lambda'$;
\item  $C_0= \mho_\Lambda$ is the set of all Delone sets in $\Omega_{\Lambda}$ containing origin.
\end{itemize}

Next, for some large $R$, consider $C_1\subset C_{0,1}$ which is the clopen subset of sets such that the $R$-cluster around the origin coincides with that of $\Lambda'$. If 
      $$\mathcal{R}_1:= \{t\in\mathbb{R}^d:\varphi_t(\hat{\Lambda})\in C_1\mbox{ for some }\hat{\Lambda}\in C_1\}$$
      is the Delone set of return vectors to $C_1$, then it defines a tiling $\mathcal{T}_1$ through its Voronoi tesselation. The crucial step at this point is to reconcile the boxes defined by the tiling $\mathcal{T}_0$ with those of $\mathcal{T}_1$. This can be done through a careful modification of the boundaries of the tiles of $\mathcal{T}_1$ (\cite[Equation (4.1)]{AC:tower}) and leads to a collection of boxes $\mathcal{B}_1$ which can be seen to be zoomed out of those defined by $C_0$, $\mathcal{B}_0$. Continuing recursively, one obtains a tower system, and the estimates of linear repetitivity give the estimates (ii) of the theorem above. This argument carries through if the starting set is a Delone multiset rather than a Delone set.
      \end{remark}
%      Observe that the bounds on the internal and external radii $r_{int}$ and $R_{ext}$ for clusters which give the $n^{th}$ box decomposition imply that for all $n>0$ and $T>\lambda^n$
%        $$ K'_3 \frac{\mathrm{Vol}(B_T(0))}{(K_2\lambda^n)^d} \leq   L^\Lambda_{B_{\lambda^n}(0)\cap \Lambda}(B_T(0)) \leq K'_4 \frac{\mathrm{Vol}(B_T(0))}{(K_1\lambda^n)^d},$$
%      where the constants $K'_3,K'_4>0$ depend only on the dimension.

      The theorem above about tower systems for linearly repetitive sets implies some ``balanced'' properties of the transversal sets $C_{n,i}$. First, we now note that since 
      \begin{equation}
        \label{eqn:one}
        \mu(C_{n,i}[D_{n,i}])= \nu(C_{n,i})\mbox{Vol}(D_{n,i})\hspace{.5in}\mbox{ and }\hspace{.5in}1 = \sum_{i=1}^{k_n} \nu(C_{n,i})\mbox{Vol}(D_{n,i}),
      \end{equation}
      picking $R = r_{int}(\mathcal{B}_n)\geq K_1\lambda^n$ and noting that $B_R(0)\subset D_{n,i}$, we have that there exists a $K_3$ such that
      $$\nu(C_{n,i}) \leq K_3 \lambda^{-dn}$$
      for all $n>0$ and $1\leq i \leq k_n$. In addition, the theorem above implies (\cite[Lemma 18]{ACCDP:LR}) the existence of a $c>0$ such that
      $$\inf_{\substack{n>0 \\ 1\leq i \leq k_n }}\mathrm{Vol}(D_{n,i})\nu(C_{n,i})>c $$
      which, when noting that $\mathrm{Vol}(D_{n,i}) < K_4' \lambda^{dn} $ for some $K_4'>1$ and all $n$, implies that
      \begin{equation}
        \label{eqn:balance}
        K^{-1}_4 \lambda^{-  dn}\leq  \nu(C_{n,i})
      \end{equation}
      for some $K_4>1$, for all $n>0$ and $1\leq i \leq k_n$. By Remark \ref{rem:general}, these estimates also hold for linearly repetitive Delone multisets.\\

Statement (ii) of the above theorem implies the following estimate, which will be crucial in proving the main result of this paper. For $n\in\mathbb{N}$ and $1\leq i  \leq k_n$  
define the set of return vectors to $C_{n,i}$ as
      $$\mathcal{R}_{n,i}:= \{t\in\mathbb{R}^d:\varphi_t(\hat{\Lambda})\in C_{n,i}\mbox{ for some }\hat{\Lambda}\in C_{n,i}\}.$$
      It is a Delone set with packing radius greater than $K_1\lambda^n-2R_\Lambda$, that is, 
      \begin{equation}
        \label{eqn:farReturns}
        K_1\lambda^n-2R_\Lambda \leq r_{\mathcal{R}_{n,i}}.
      \end{equation}

      We now introduce a condition which will simplify things significantly.
      \begin{hypo}
        \label{hyp:1}
        The nice tower system for a linearly repetitive Delone multiset has $k_n>1$ for infinitely many $n$.
      \end{hypo}
      \begin{remark}
        Hypothesis \ref{hyp:1} holds in the linearly repetitive examples that are prominent in the literature, namely self-similar tilings/sets and globally random substitution tilings/sets (e.g. \cite{ST:random}). More specifically, for substitution systems, $k_n$ can be constant and corresponds to the number of prototiles in the tiling. This is also true for the globally random substitution systems in \cite{ST:random}. Thus, asking that $k_n>1$ for infinitely many $n$ in Hypothesis \ref{hyp:1} is not asking for too much.
      \end{remark}

      \section{Far away functions}
      \label{sec:Far}
      Consider an aperiodic, linearly repetitive collection of scatterers $\mathcal{S}\subset\mathbb{R}^d$. Without loss of generality we assume that the center of mass of one of the scatterers coinsides with the origin, and denote their corresponding aperiodic, linearly repetitive Delone multiset $\Lambda_\mathcal{S}$, its nice tower system $\mathcal{B}_\mathcal{S}$.
      Using the coordinates (\ref{eqn:coordinates}), for a function $\psi:\mathcal{D}\rightarrow \mathbb{R}$, analogous to (\ref{eqn:HoldConst}), let
      \begin{equation}
        \label{eqn:HoldConst2}
        |\psi|_\alpha^\perp:= \max_{i\in\mathfrak{C}}\sup_{(x,v)\in T^1_+(\partial S_i)}\sup_{\substack{c_1\neq c_2 \\ \in \mathcal{C}_i}} \frac{\left| \psi(x,v,c_1)-\psi(x,v,c_2)\right|}{d(c_1,c_2)^\alpha},
      \end{equation}
      where $d(c_1,c_2)$ is the metric inherited from (\ref{eqn:distance}) restricted to local transversals.
      Let $H_\alpha^\perp(\mathcal{S})$ be the space of continuous functions $\psi:\mathcal{D}\rightarrow \mathbb{R}$ which are transversally $\alpha$-H\"older functions, that is, which have $|\psi|^\perp_\alpha<\infty$. Endowed with the norm $\|\psi\|_\alpha = \|\psi\|_{C^0}+|\psi|^\perp_\alpha$, it is a Banach space.
      
      Using the coordinates (\ref{eqn:coordinates}), define the sets
      $$A_{n,i} := T^1_+(\partial S_{j_{n,i}})\times  C_{n,i}\subset \mathcal{D},$$
      where $j_{n,i}\in\mathfrak{C}$ is the index associated to the point defining the base of the transversal in $C_{n,i}[D_{n,i}]$. Here and in what follows we only take $n$ such that $k_n>1$ as in Hypothesis \ref{hyp:1}. Define the functions
      $$\psi_{i}^{(n)} := \chi_{A_{n,i}}.$$
      We want to estimate the transversal H\"older constant (\ref{eqn:HoldConst2}) of these functions. %, that is, we are interested in bounding
%      $$|\psi_i^{(n)}|_{\alpha}^\perp:= \sup_{\Lambda'\in\Omega_\Lambda}\sup_{\substack{\Lambda_1\neq\Lambda_2 \\ \in C_{\Lambda',R_\Lambda}}} \frac{|\psi^{(n)}_i(\Lambda_1) - \psi^{(n)}_i(\Lambda_2)|}{d(\Lambda_1,\Lambda_2)^\alpha}.$$
      %for $\alpha>0$.
      By definition, the quantity $|\psi^{(n)}_i(p_1) - \psi^{(n)}_i(p_2)|/d(p_1,p_2)^\alpha$ is nonzero only when $p_1$ is in $A_{n,i}$ and $p_2$ is not (or vice-versa). Thus we are interested in knowing how transversally close $p_1$ and $p_2$ can be while having one in $A_{n,i}$ and the other not.

      Let $p_1\neq p_2\in\mathcal{D}$ have the same $T^1_+(\partial\mathcal{S}_i)$ coordinates, so that the distance between them is measured entirely in the transverse direction as they are in the same local transversal. By the estimates on the box decompositions of the tower system, if they have the same $(2L+1)\lambda^n$ neighborhood around the origin, then they would have to be in the same box. Thus if $p_1\in A_{n,i}$ and $p_2\not\in A_{n,i}$ (or vice-versa), then $d(p_1,p_2)\geq ((2L+1)\lambda^n)^{-1}$. Thus it follows that for any $i$ we have
      \begin{equation}
        \label{eqn:HolderEst}
        |\psi_i^{(n)}|^\perp_\alpha \leq (2L+1)^\alpha \lambda^{\alpha n}.
      \end{equation}
      
      \section{Non-fast mixing}
      \label{sec:nonFast}
      %Denote $\Lambda=\Lambda_\mathcal{S}$. For any cluster of the form $P = B_R(p)\cap \Lambda$, define $\Lambda_P$ to be the Delone set consisting of points which are centers of clusters which are translation equivalent to $P$. Let $r_P,R_P$ be the packing and covering radii, respectively, of $\Lambda_P$. 
      
      %Let $P_1 = B_R(p_1)\cap \Lambda$ be a large cluster for some $p_1\in\Lambda$. It follows from Lemma \ref{lem:rR} that for any $p_2\in B_{\frac{R}{3(L+1)}}(p_1)\cap \Lambda$, the cluster $P_2 = B_R(p_2)$ is not translation equivalent to $P_1$. For $i = 1,2$, using the coordinates (\ref{eqn:coordinates}), define the sets

      Here we prove the main result. Let $\psi_i^{(n)}$ and $\psi_j^{(n)}$ be two functions constructed from a nice box decomposition as in the previous section and consider the correlation
      $$\int_\mathcal{D} \psi_{i}^{(n)}\circ F^k \cdot \psi_{j}^{(n)}\, d\mu,$$
      where $\mu$ is an $F$-invariant probability measure on $\mathcal{D}$, which is absolutely continuous in the sense of Definition \ref{def:abscont}.
%locally (using coordinates (\ref{eqn:coordinates})) of the form $L_i\times \nu$, where $L_i$ is an absolutely continuous measure on $T^+_1(\partial S_i)$, and $\nu$ is the unique frequency measure defined in (\ref{eqn:freqMeas}).
      Note that by Hypothesis \ref{hyp:1} it can be assumed that $i\neq j$, so let us assume $i\neq j$ for the rest of the section. Let $P_{n,i}$ be the clusters defined by the tower system, that is,
      $$ P_{n,i} := \{t\in D_{n,i}:\varphi_t(\Lambda')\in \mho_\Lambda\mbox{ for every $\Lambda'\in C_{n,i}$}\}. $$
      Since the function $\psi_{i}^{(n)}$ is supported on scatterers whose $(K_1\lambda^n- 2R_\Lambda)$-neighborhood is determined by the cluster $P_{n,i}$, and the distance between two points at the centers of different clusters is at least $2(K_1\lambda^n-2R_\Lambda)$ by (\ref{eqn:farReturns}), $\psi_{i}^{(n)}\circ F^k(x) \cdot \psi_{j}^{(n)}(x) = 0$ for all the $k$ so that no orbit starting at the center of a cluster of the form $P_{n,j}$ reaches a scatterer at the center of the cluster of the form $P_{n,i}$ after $k$ collisions. By finite horizon, a trajectory can travel at most a distance of $Mk$ after $k$ collisions. Thus, by (\ref{eqn:farReturns}),
      $$\psi_{i}^{(n)}\circ F^k(x) \cdot \psi_{j}^{(n)}(x) = 0 \hspace{.45in}\mbox{ for all } \hspace{.45in}|k|< \frac{2(K_1\lambda^n-2 R_\Lambda- B_S)}{M},$$
      where $B_S$ is any positive constant such that for any $i\in\mathfrak{C}$ and $x\in\partial S_i$, the distance between $x$ and the center of mass of $S_i$ is less than or equal to $B_S$. The constant exists since $\mathfrak{C}$ is finite. Thus we have that
      \begin{equation}
        \label{eqn:bnd1}
        \left|\int_\mathcal{D} \psi_{i}^{(n)}\circ F^k \cdot \psi_{j}^{(n)} \, d\mu - \mu(\psi_i^{(n)})\mu(\psi_j^{(n)})\right| = \left| \mu(\psi_i^{(n)})\mu(\psi_j^{(n)})\right| \geq \varrho^2 K_4^{-2}\lambda^{-2dn}
      \end{equation}
      for $|k|<  \frac{2}{M}(K_1\lambda^n - 2R_\Lambda - B_S)$ by (\ref{eqn:balance}), where $\varrho = \min_i L_i(T^1_+(\partial S_i))$.

      We now show how exponential mixing is incompatible with the nice tower system which linearly repetitive systems have. Suppose that for some $\alpha>0$ and $\tau\in (0,1)$ there is a $C>0$ such that
      $$\left|\int_\mathcal{D} h_1\circ F^k \cdot h_2 \, d\mu - \mu(h_1)\mu(h_2)\right|\leq C \|h_1\|_\alpha \|h_2\|_\alpha\tau^k$$
      for all $k\in\mathbb{N}$ and $h_i\in H_\alpha$. Using the functions $\psi_i^{(n)}$ we constructed above, by (\ref{eqn:HolderEst}) and (\ref{eqn:bnd1}), we would have
      $$\varrho^2 K_4^{-2}\lambda^{-2dn} \leq C\|\psi_{i}^{(n)}\|_\alpha\|\psi_{j}^{(n)}\|_\alpha\tau^k \leq  C (1+(2L+1)^{\alpha} \lambda^{\alpha n})^2\tau^k$$
      for all $|k|<  \frac{2}{M}(K_1\lambda^n - 2R_\Lambda - B_S)$. In particular this should hold for $k=\left\lfloor \frac{2}{M}(K_1\lambda^n - 2R_\Lambda- B_S)\right\rfloor$, so
      $$\varrho^2 K_4^{-2}\lambda^{-2dn} \leq C (1+(2L+1)^\alpha \lambda^{\alpha n})^2\tau^{\left\lfloor\frac{2}{M}(K_1\lambda^n - 2R_\Lambda)\right\rfloor}.$$
      This implies
      $$\log(C^*)-2n(d+\alpha)\log\lambda\leq \left\lfloor\frac{2}{M}(K_1\lambda^n-2R_\Lambda)\right\rfloor\log \tau$$
      for some $C^*>0$.  Now, the right hand side is negative and decreasing at an exponential rate with $n$, while the left hand side is negative and decreasing at a linear rate. This of course cannot happen for all $n$, and so this obstructs the system from mixing at an exponential rate. The same argument prevents the system from mixing with any stretched exponential rate.

      Polynomial mixing is equally restricted by the nice tower system which linearly repetitive systems have. Suppose that for some $\alpha>0$ and $\gamma>0$ there is a $C>0$ such that
      $$\left|\int_\mathcal{D} h_1\circ F^k \cdot h_2 \, d\mu - \mu(h_1)\mu(h_2)\right|\leq C \|h_1\|_\alpha \|h_2\|_\alpha |k|^{-\gamma}$$
      for all $k\in\mathbb{Z}$ and $h_i\in H_\alpha$. Using the functions $\psi_i^{(n)}$ we constructed above, by (\ref{eqn:HolderEst}) and (\ref{eqn:bnd1}), we would have
      $$\varrho^2 K_4^{-2}\lambda^{-2dn} \leq C\|\psi_{i}^{(n)}\|_\alpha\|\psi_{j}^{(n)}\|_\alpha k^{-\gamma} \leq  C (1+(2L+1)^\alpha \lambda^{\alpha n})^2k^{-\gamma}$$
      for all $|k|<  \frac{2}{M}(K_1\lambda^n - 2R_\Lambda- B_S)$. In particular this should hold for $k=\left\lfloor \frac{2}{M}(K_1\lambda^n - 2R_\Lambda-B_S)\right\rfloor$, so
      $$\varrho^2 K_4^{-2}\lambda^{-2dn} \leq C(1+ (2L+1)^\alpha\lambda^{\alpha n})^2\left\lfloor\frac{2}{M}(K_1\lambda^n - 2R_\Lambda- B_S)\right\rfloor^{-\gamma}.$$
      For sufficiently large $n$ this implies the inequality
      $$C^{**} \lambda^{-2(d+\alpha)n}\leq \lambda^{-\gamma n}.$$
      This cannot hold for all $n$ unless $\gamma\leq 2(d+\alpha)$ . 
      \bibliographystyle{amsalpha}
      \bibliography{biblio}
\end{document}